\documentclass[11pt,fleqn]{article}
\usepackage{amssymb}
\usepackage{amsfonts}
\usepackage{amsmath}
\usepackage{amsthm}
\usepackage{graphicx}
\usepackage{epsfig}
\usepackage{psfrag}
\usepackage{color}
\usepackage{dsfont}
\makeatletter
\xdef\@endgadget#1{{\unskip\nobreak\hfil\penalty50\hskip1em\hbox{}\nobreak
    \hfil#1\parfillskip=0pt\finalhyphendemerits=0\par}}
\def\@qedsymbol{${}_\blacksquare$}
\def\qed{\@endgadget{\@qedsymbol}}
\newtheorem{lemma}{Lemma}[section]
\newtheorem{theorem}[lemma]{Theorem}

\newtheorem{assumption}[lemma]{Assumption}
\newtheorem{proposition}[lemma]{Proposition}
\newtheorem{remark}[lemma]{Remark}
\newcommand{\mR}{\mathbb{R}}
\newcommand{\mC}{\mathbb{C}}

\newcommand{\bq}{\begin{equation}}
\newcommand{\eq}{\end{equation}}

\newcommand{\D}{\mathcal{D}}
\renewcommand{\L}{\mathcal{L}}

\DeclareMathOperator{\spa}{span}

\def\BibTeX{{\rm B\kern-.05em{\sc i\kern-.025em b}\kern-.08em
    T\kern-.1667em\lower.7ex\hbox{E}\kern-.125emX}}

\title{\LARGE \bf The flow equations of linear resistive electrical networks}

\author{Arjan van der Schaft
\thanks{A.J. van der Schaft is with the Johann Bernoulli Institute for Mathematics and Computer
Science, and the Jan C. Willems Center for Systems and Control, University of Groningen, PO Box 407, 9700 AK, the
Netherlands,
        {\tt\small A.J.van.der.Schaft@rug.nl}}
}
\date{}

\begin{document}

\maketitle
\thispagestyle{empty}
\pagestyle{empty}


\section{Introduction}
The theory of electrical engineering has triggered and stimulated many interesting and elegant mathematical developments, which conversely turned out to be of great importance for their engineering applications. Especially within systems and control theory we can witness a wide range of examples of this fruitful interplay; from network synthesis, impedance matching, uncertainty reduction, robust control design, to large-scale network systems.

The current paper focusses on a particularly classical, and in some sense very simple, example; namely the study of {\it linear resistive electrical networks}. The mathematical formulation of resistive networks, and the study of their properties, goes back to the 1847 paper \cite{Kirchhoff} (while still being a student) by Gustav Robert Kirchhoff (1824-1887), which laid the foundations for the theory of electrical networks, and at the same time can be seen as one of the early milestones in the development of (algebraic) graph theory. From a systems theory point of view it is a beautiful example of an {\it interconnected} and {\it open} system; although there is no dynamics associated with it. See also \cite{willemsportsterminals} for stimulating reading, as well as \cite{NOW} for related developments\footnote{While writing the present paper I came across the recent paper \cite{doerfler1}, surveying, rather complementarily to the present paper, the interplay between algebraic graph theory and electrical networks from many angles.}.

After recalling in Section 2 some basic notions of algebraic graph theory, Section 3 will summarize the graph-theoretic formulation of resistive networks. Section 4 starts with the description and solution of the classical problem posed by Kirchhoff concerning the existence and uniqueness of solutions of the flow equations of open resistive networks, and discusses some related problems. Although it is very difficult to claim any originality in this very classical area (with the names of e.g. Maxwell, Rayleigh, Thomson and Dirichlet attached to it), I believe there is still some originality in Propositions \ref{prop1} and \ref{prop2}.

Section 5 discusses within the same framework the reconstructibility of the conductance parameters of a resistive network from the knowledge of the map from boundary nodal voltages potentials to boundary nodal currents. This problem is a discrete version of the well-known inverse boundary problem of Calderon \cite{calderon}.
Finally, Section 6 indicates extensions to RLC-networks, and in particular the very important power-flow problem in electrical power networks. Section 7 contains the conclusions. 

%

\section{Preliminaries about graphs}
First, let us recall from e.g. \cite{bollobas, godsil} the following standard notions regarding directed graphs. 
A {\it graph} is defined by a set of {\it nodes} (or vertices), and a set of {\it edges} (or links, branches), where the set of edges is identified with a subset of unordered pairs of nodes; $\{i,j\}$ denotes an edge between nodes $i$ and $j$. We allow for multiple edges between nodes, but not for self-loops $\{i,i \}$. By endowing the edges with an orientation we obtain a {\it directed graph}, with edges corresponding to ordered pairs; $ij$ denotes an edge from node $i$ to node $j$. In the sequel, 'graph' will always mean 'directed graph'.

A directed graph with $N$ nodes and $M$ edges is specified by its $N \times M$ {\it incidence matrix}, denoted by $\D$. Every column of $\D$ corresponds to an edge of the graph, and contains exactly one $1$ at the row corresponding to the tail node of the edge and one $-1$ at the row corresponding to its head node, while the other elements are $0$. In particular, $\mathds{1}^T\D=0$ where $\mathds{1}$ is the vector of all ones. 
Furthermore, $\ker \D^T = \spa \mathds{1}$ if and only if the graph is {\it connected} (any node can be reached from any other vertex by a sequence of, - undirected -, edges). In general, the dimension of $\ker \D^T$ is equal to the number of connected components. Throughout this paper we will adopt\footnote{Without loss of generality since otherwise the analysis can be done for every connected component of the graph.}
\begin{assumption}
The graph is connected.
\end{assumption}

Corresponding to a directed graph we can define, in the spirit of general $k$-complexes, the following vector spaces; see \cite{hamgraphs}. The {\it node space} $\Lambda_0$ of the directed graph is defined as the set of all functions from the node set $\{1,\cdots,N\}$ to $\mR$. Obviously $\Lambda_0$ can be identified with $\mR^N$. The dual space of $\Lambda_0$ is denoted by $\Lambda^0$. Furthermore, the {\it edge space} $\Lambda_1$ is defined as the linear space of functions from the edge set $\{1,\cdots,M \}$ to $\mR$, with dual space denoted by $\Lambda^1$. Both spaces can be identified with $\mR^M$. It follows that the incidence matrix $\D$ defines a linear map (denoted by the same symbol) $\D: \Lambda_1 \to \Lambda_0$ with adjoint map $\D^T: \Lambda^0 \to \Lambda^1$. In the context of the present paper the space $\Lambda_1$ corresponds to the {\it currents} through the edges, and the dual space $\Lambda^1$ to the space of {\it voltages} across the edges. Furthermore, as we will see later on, $\Lambda_0$ denotes the space of {\it nodal currents} (entering the nodes from the environment, or from electrical devices located at the nodes), and $\Lambda^0$ the space of {\it voltage potentials} at the nodes.

Finally, we note that it is straightforward to extend the network models described in this paper, as well as their dynamical versions (see e.g. \cite{hamgraphs}), to more general cases. Indeed, for any linear space $\mathcal{R}$ (e.g., $\mathcal{R} = \mR^3$) we can instead define $\Lambda_0$ as the set of functions from $\{1,\cdots,N\}$ to $\mathcal{R}$, and $\Lambda_1$ as the set of functions from $\{1,\cdots,M \}$ to $\mathcal{R}$. In this case we identify $\Lambda_0$ with the tensor product $\mR^N \otimes  \mathcal{R}$, and $\Lambda_1$ with the tensor product $\mR^M \otimes  \mathcal{R}$. Furthermore, the incidence matrix $\D$ defines a linear map $\D \otimes I : \Lambda_1 \to \Lambda_0$, where $I$ is the identity map on $\mathcal{R}$. In matrix notation $\D \otimes I$ equals the Kronecker product of the incidence matrix $\D$ and the identity matrix $I$. See \cite{bosgra, hamgraphs} for further details. 

\section{Resistive electrical networks}
Consider a linear resistive electrical network, where the underlying connected graph, with $N$ nodes and $M$ edges, is specified by its incidence matrix $\D$. 
Then the relation between the vector of {\it voltage potentials} 
\bq
\psi=(\psi_1, \cdots, \psi_N)^T
\eq
at the nodes, and the vector of {\it nodal currents}
\bq
J=(J_1,\cdots,J_N)^T
\eq
entering the nodes (e.g., the currents taken from generators or delivered to loads) is determined as follows. By Kirchhoff's {\it voltage laws}, the vector of {\it voltages} $V=(V_1,\cdots,V_M)^T$ across the edges is determined by the vector of nodal voltage potentials as
\bq
V=\D^T\psi
\eq
Dually, the relation between the vector of nodal currents $J$ and the vector $I=(I_1, \cdots,I_N)^T$ of currents through the edges is given by Kirchhoff's {\it current laws}
\bq
J= \D I
\eq
It follows that $I \in \Lambda_1$, $V \in \Lambda^1$, $\psi \in \Lambda^0$, and $J \in \Lambda_0$. In particular, the products $V^TI$ and $\psi^TJ$ are intrinsically defined, while
\bq
\label{power}
V^TI=\psi^T \D I = \psi^TJ,
\eq
expressing that the total power through the edges is equal to the total power at the nodes.

\smallskip

In the case of a linear resistive network each edge corresponds to a {\it linear resistor}. That is, for the $k$-th edge the current $I_k$ through the edge, and the voltage $V_k$ across the edge, are related as
\bq
I_k= g_kV_k ,
\eq
with $g_k$ the {\it conductance} of the resistor, i.e., $g_k=\frac{1}{r_k}$ with $r_k>0$ the resistance of the $k$-th edge. Hence 
\bq
I=GV
\eq
with $G$ the $M \times M$ diagonal matrix of conductances. Putting this all together, one obtains
\bq
J= \D G \D^T\psi
\eq
The matrix $\L=\D G \D^T$ is called the {\it Laplacian matrix}\footnote{Classically, in the context of resistive networks the matrix $\L$ is called the {\it Kirchhoff matrix}; however we will stick to the currently used terminology in general network models.} of the electrical network. Laplacian matrices are fully characterized by the following properties.
\begin{proposition}\cite{Schaft_SCL10}
\label{prop:laplacian}
The Laplacian matrix $\L$ is a symmetric matrix, with all its diagonal elements {\it positive}, and all off-diagonal elements non-positive. Furthermore, the row and column sums of $\L$ are zero, i.e., $\L \mathds{1}=0, \mathds{1}^T\L=0$. 
\end{proposition}
'Fully characterized' means that conversely for any matrix $\L$ satisfying the properties listed in Proposition \ref{prop:laplacian} there exists an incidence matrix $\D$ and positive diagonal matrix $G$ such that $\L=\D G \D^T$. In particular, any non-zero off-diagonal element $\L_{ij}$ defines an edge $ij$ with conductance equal to $-\L_{ij}$.
Furthermore, it is well-known (see e.g. \cite{bollobas}) that $\L$ is independent of the orientation of the graph. Thus if we take another orientation of the graph (corresponding to an incidence matrix $\D$ where some of the columns have been multiplied by $-1$), then $\L$ remains the same. Moreover, by the assumption of connectedness $\ker \L=\spa \mathds{1}$.

\smallskip

Also note that any Laplacian matrix $\L$ is positive-semidefinite; in fact, the quadratic form associated to $\L$ is equal to
\bq
\psi^T\L \psi = \sum_{ij}g_{ij}(\psi_i-\psi_j)^2,
\eq
where the summation is over all the $M$ edges $ij$. Note that $\psi^T\L \psi=\psi^TJ$ is equal to the total power at the nodes, which by \eqref{power} is equal to the total power $V^TI$ at the edges.

\begin{remark}
The same equations hold for other classes of physical systems \cite{sch-willems}. For example, mechanical damper systems are completely analogous, with the voltage potentials $\psi$ replaced by the {\it velocities} of the nodes, the nodal currents $J$ by the nodal {\it forces}, and the $k$-the edge corresponding to a linear damper $F_k=d_kv_k$, with $d_k$ the damping constant, $F_k$ the damping force, and $v_k$ the difference of the velocities of the tail and head node of this edge.
\end{remark}

\section{Kirchhoff's problem and its dual}
Now suppose the $N$ nodes of the network graph are {\it split} into $N_C$ internal {\it connection} nodes (denoted by $C$), and $N_B$ remaining external {\it boundary} nodes (denoted by $B$), with $N=N_B + N_C$. Correspondingly we split (possibly after reordering) the vectors $J$ and $\psi$ of nodal currents and voltage potentials as
\bq
J= \begin{bmatrix}J_B \\ J_C \end{bmatrix}, \quad \psi= \begin{bmatrix}\psi_B\\ \psi_C \end{bmatrix}
\eq
and the incidence matrix into
\bq
\D=\begin{bmatrix} \D_B\\ \D_C \end{bmatrix}
\eq
One then obtains the following equations
\bq
\label{1}
\begin{bmatrix} J_B \\ J_C \end{bmatrix} = 
\L
\begin{bmatrix} \psi_B \\ \psi_C \end{bmatrix},
\eq
where
\bq
\L= \begin{bmatrix} \D_B \\ \D_C \end{bmatrix} G \begin{bmatrix} \D^T_B & \D^T_C \end{bmatrix} =:
\begin{bmatrix} \L_{BB} & \L_{BC} \\ \L_{CB} & \L_{CC} \end{bmatrix}
\eq

\subsection{Kirchhoff's problem}
The classical {\it Kirchhoff problem} and its solution, dating back to \cite{Kirchhoff}, can be formulated in this notation as follows. Fix the boundary voltage potentials $\psi_B=\psi_B^* \in \mR^{N_B}$, and let $0=J_C \in \mR^{N_C}$ (no nodal currents at the connection nodes). Then there exists \cite{bollobas} a {\it unique} vector of voltage potentials $\psi_C^* \in \mR^{N_C}$ at the connection nodes such that
\bq
\begin{bmatrix} J^*_B \\ 0 \end{bmatrix} = 
\begin{bmatrix} \L_{BB} & \L_{BC} \\ \L_{CB} & \L_{CC} \end{bmatrix}
\begin{bmatrix} \psi^*_B \\ \psi^*_L \end{bmatrix} ,
\eq
with the corresponding vector of currents $J_B^* \in \mR^{N_B}$ at the boundary nodes given by\footnote{It can be easily seen from the assumption of connectedness that the submatrix $\L_{CC}$ is invertible.}
\bq
J_B^*= (\L_{BB} - \L_{BC}\L_{CC}^{-1}\L_{CB})\psi_B^*=\L_S \psi_B^*,
\eq
with $\L_S:=(\L_{BB} - \L_{BC}\L_{CC}^{-1}\L_{CB})$ the Schur complement of the Laplacian matrix $\L$ with respect to $\L_{CC}$. 

A key observation is that the Schur complement $\L_S$ is again a Laplacian matrix. Explicit statements of this result in the literature before its formulation and proof in \cite{Schaft_SCL10} seem hard to find; although \cite{crab} contains the closely related result that the Schur complement of any $M$-matrix is again an $M$-matrix\footnote{I thank Nima Monshizadeh for pointing out this reference to me.}. See also \cite{rommes} for the use of Schur complements in reduction of large-scale resistive networks.
\begin{theorem} \cite{Schaft_SCL10}
If the graph $\mathcal{G}$ is connected, then all diagonal elements of $\L=\D G \D^T$ are $>0$. Furthermore, all Schur complements of $\L=\D G \D^T$ are well-defined, symmetric, with diagonal elements $>0$, off-diagonal elements $\leq 0$, and with zero row and column sums. In particular, all Schur complements $\L_S$ of $\L=\D G \D^T$ are Laplacian, and thus can be written as $\bar{\D} \bar{G} \bar{\D}^T$, with $\bar{\D}$ the incidence matrix of some connected graph with the same set of nodes, and $\bar{G}$ a positive definite diagonal matrix with the same dimension as the number of edges of the graph defined by $\bar{\D}$. 
\end{theorem}
The proof given in \cite{Schaft_SCL10} is based on two observations. First, the Schur complement of a Laplacian matrix $\L$ with respect to a scalar diagonal element can be immediately checked to be again a Laplacian matrix; see also \cite{willemsverriest1}. Second, any Schur complement can be obtained by the successive application of taking Schur complements with respect to diagonal elements. In fact, this follows from the quotient formula given in \cite{crab}. Denote by $M/P$ the Schur complement of the square matrix $M$ with respect to a leading square submatrix $P$. Then for any leading submatrix $Q$ of $P$ we have the equality
\bq
M/P = (M/Q)/(P/Q)
\eq
Thus $\L_S$ is again a Laplacian matrix. Therefore we have obtained a {\it reduced} resistive network with only boundary nodes, where the map from boundary nodal voltages to boundary nodal currents is given by
\bq
J_B = \L_S \psi_B
\eq
This reduced resistive network is {\it equivalent} (as seen from the boundary nodes) to the original one with the constraint $J_C=0$. The transformation from electrical network with boundary and internal nodes satisfying $J_C=0$ to the reduced network {\it without} internal nodes is called {\it Kron reduction} \cite{kron}; see e.g. \cite{doerfler} for a review of the literature and applications of Kron reduction.

\begin{remark}
A {\it special} case is obtained by considering just {\it two} boundary nodes; the rest being connection nodes. Then the $2 \times 2$ Schur complement $\L_S$ is of the form
\bq
\L_S= \begin{bmatrix} g &- g \\-g & g \end{bmatrix}
\eq
with $g$ the {\it effective conductance} between the two boundary nodes, and $R=\frac{1}{g}$ the {\it effective resistance}.
\end{remark}

The uniquely determined vector of internal voltage potentials $\psi_C^*$, called the {\it open-circuit} internal voltage potentials, corresponding to the boundary voltage potentials $\psi_B^*$, is given as
\bq
\label{kirchhoffmap}
\psi_C^*= - \L_{CC}^{-1}\L_{CB}\psi_B^*
\eq
Since $\L_{CC}$ is an invertible $M$-matrix it follows that $\L_{CC}^{-1}$ has all nonnegative elements, and thus also $- \L_{CC}^{-1}\L_{CB}$ is a matrix with all nonnegative elements. Furthermore, it follows from the discrete {\it Maximum Modulus principle} \cite{bollobas} that for all $\psi_B^*$
\bq
\|\psi_C^*\|_{\mathrm{max}} \leq \|\psi_B^*\|_{\mathrm{max}},
\eq
where $\| \cdot \|_{\mathrm{max}}$ denotes the max-norm: $\|\psi\|_{\mathrm{max}}=\max_j |\psi_j|$. Hence the linear map $\psi_B^* \mapsto \psi_C^*$ given by the matrix $- \L_{CC}^{-1}\L_{CB}$ has induced norm $\leq 1$. In fact, since $\psi_B^* = \mathds{1}$ yields $\psi_C^*= \mathds{1}$, it follows that
\bq
\|- \L_{CC}^{-1}\L_{CB}\|_{\mathrm{max}} =1
\eq
The {\it open-circuit} internal voltage potential $\psi_C^*$ has the following classical {\it minimization} interpretation, called Thomson's (or also Dirichlet's) principle \cite{bollobas, Schaft_SCL10}. Given $\psi_B^*$, then $\psi_C^*$ is the unique minimizer of
\bq
\min_{\psi_C} \begin{bmatrix} \psi^{*T}_B & \psi^T_C \end{bmatrix} \L\begin{bmatrix} \psi^*_B \\ \psi_C \end{bmatrix} = \min_{\psi_C} \sum_{(ij)} g_{(ij)}(\psi_i-\psi_j)^2 
\eq
(with $g_{(ij)}$ denoting the conductance of the edge $ij$, and summation over all the edges $ij$ of the circuit graph).
Indeed, the gradient vector of this expression with respect to $\psi_C$ is equal to two times $I_C= \L_{CB}\psi_B^*+\L_{CC}\psi_C$. The expression $\sum g_{(ij)}(\psi_i-\psi_j)^2$ equals the dissipated power in the resistive network. Thus the open-circuit internal voltage potential $\psi_C^*$ corresponds to {\it minimal dissipated power} in the resistive network.

This is summarized in the following proposition.
\begin{proposition}
\label{prop1}
The linear map $\psi_B^* \mapsto \psi_C^*$ given by \eqref{kirchhoffmap} has induced $\mathrm{max}$-norm equal to $1$. Given $\psi_B^*$, then $\psi_C^*$ is the unique minimizer of the dissipated power $\sum g_{(ij)}(\psi_i-\psi_j)^2$. Any other resistive network inducing by Kron reduction the same linear map $J_B = \L_S \psi_B$ has the same dissipated power.
\end{proposition}
\subsection{The dual problem}
Dually to Kirchhoff's problem there is the (easier) {\it short-circuit} problem, corresponding to taking the boundary voltage potentials $\psi_B=0$, and instead to consider a prescribed value of the vector of internal nodal currents $J_C=\bar{J}_C$. Then the corresponding internal voltage potentials $\bar{\psi}_{C}$ are given as
\bq
\label{kirchhoffmap1}
\bar{\psi}_{C} = \L_{CC}^{-1}\bar{J}_C ,
\eq
while the short-circuit boundary nodal currents $\bar{J}_B$ are determined by $\bar{J}_C$ as
\bq
\bar{J}_B= \L_{BC}\L_{CC}^{-1} \bar{J}_C
\eq
Note that the matrix $\L_{BC}\L_{CC}^{-1}$ has all non-positive elements, and equals minus the transpose of the matrix $\L_{CC}^{-1}\L_{CB}$ of the Kirchhoff problem. Hence the map $\bar{J}_C \mapsto - \bar{J}_B$ is {\it dual} to the map $\psi_B^* \mapsto \psi_C^*$. Thus\footnote{\eqref{dual} also follows from Tellegen's theorem.}
\bq
\label{dual}
\bar{J}^T_C \psi_C^*= - \bar{J}^T_B \psi_B^*
\eq
for any $\psi_B^*,\bar{J}_C$ and corresponding $\psi_C^*,\bar{J}_B$.
In particular the induced $1$-norm of $\L_{BC}\L_{CC}^{-1}$ is equal to $1$:
\bq
\|\L_{BC}\L_{CC}^{-1}\|_1=1
\eq
This follows from considering $\psi_B^*=\mathds{1}$ and $\psi_C^*=\mathds{1}$, and computing
\bq
\|\bar{J}_B \|_1= |\bar{J}^T_B \mathds{1}| =|\bar{J}^T_B \psi_B^*|=|\bar{J}^T_C \psi_C^*|=\|\bar{J}_C \|_1
\eq
for all $\bar{J}_C$ and corresponding $\bar{J}_B$.
This is summarized in the following proposition.
\begin{proposition}
\label{prop2}
The map $\bar{J}_C \mapsto - \bar{J}_B$ given by \eqref{kirchhoffmap1} has induced $1$-norm equal to $1$. It is dual to the map $\psi_B^* \mapsto \psi_C^*$ given by \eqref{kirchhoffmap}.
\end{proposition}

\subsection{Combining Kirchhoff's problem and its dual}
We can {\it combine} Kirchhoff's problem and its dual. Prescribe $\psi_B^*$ and $\bar{J}_C$. Then 
\bq
(\psi_B^*,\psi_C^*+\bar{\psi}_C, J_B^*+\bar{J}_B,\bar{J}_C)
\eq
with $\psi_C^*$ given by \eqref{kirchhoffmap}, and $\bar{J}_B$ given by \eqref{kirchhoffmap1}, is the unique solution of the resulting flow equation. 
Furthermore, this solution is the minimizer of
\bq
\min_{\psi_C} \begin{bmatrix} \psi^{*T}_B & \psi^T_C \end{bmatrix} \L\begin{bmatrix} \psi^*_B \\ \psi_C \end{bmatrix} - 2\psi_C^T\bar{J}_C
\eq
This last fact follows by noting that the gradient vector (with respect to $\psi_C$) of this expression is given by $2(J_C-\bar{J}_C)$.


\subsection{The prescribed power problem}
Still another version of the flow problem for linear resistive networks is the following. Next to the voltage potentials $\psi_B^*$ at the boundary nodes prescribe in this case the {\it power} $P_j=\psi_jJ_j$ at each of the connection nodes $j=N_B+1,\cdots, N_B+N_C=N$. Defining for every vector $z \in \mR^n$ the $n \times n$ matrix $[z]$ as the diagonal matrix with diagonal elements $z_1, \cdots,z_n$, the problem with prescribed boundary node voltages and prescribed connection nodal powers can be succinctly formulated as finding $\psi_C$ satisfying
\bq
\label{powerflow}
\bar{P}_C=[\psi_C]J_C= [\psi_C](\L_{CC}\psi_C +\L_{CB}\psi^*_B)= [\psi_C]\L_{CC}(\psi_C-\psi_C^*)
\eq
Being quadratic, this equation has generally multiple solutions.

In terms of the deviation $\psi_C-\psi_C^*$ with respect to the open-circuit load voltage potentials $\psi_C^*$ equation \eqref{powerflow} can be further rewritten as
\bq
\label{pf}
[\psi_C-\psi_C^*] \L_{CC}(\psi_C-\psi_C^*) +[\psi_C^*] \L_{CC}(\psi_C-\psi_C^*)-\bar{P}_C=0
\eq
Furthermore, we can {\it scale}, see also \cite{simpson1}, the equation \eqref{pf} by defining the column vector $x \in \mR^{N_C}$ with $i$-th element given as
\bq
x_j:=\frac{\psi_j - \psi^*_j}{\psi^*_j}, \quad j=N_B+1,\cdots, N_B+N_C=N
\eq
Then equation \eqref{pf} corresponds to
\bq
\left( [x] \L_{CC}[x] - \L_{CC}[x] + [\psi_C^*]^{-1}[\bar{P}_C] [\psi_B^*]^{-1} \right)\mathds{1}=0
\eq
in the unknown diagonal matrix $[x]$. This equation bears obvious similarities with a matrix {\it Riccati equation}.

The {\it complex} version of this problem, of much interest in the theory of {\it power networks}, will be indicated in Section 6.

\section{The discrete version of Calderon's inverse problem}

The map $J_B=\L_S\psi_B$ from boundary nodal voltage potentials to boundary nodal currents can be regarded as the {\it discrete analog} of the Dirichlet-to-Neumann map used e.g. in tomography. This map is described by the partial differential equations
\bq
\nabla \cdot (\gamma \nabla u)=0 \: \mbox{ in } \Omega,
\eq
in the unknown function $u:\Omega \to \mR$, with $\Omega$ some bounded domain in $\mR^n$, together with the Dirichlet boundary conditions 
\bq
u\mid _{\partial \Omega} =v,
\eq
where $v$ is a prescribed voltage potential function on the boundary $\partial \Omega$ of the domain $\Omega$. 
Here the function $\gamma: \Omega \to \mR$ denotes a conductance function. The Neumann boundary variables $j:= \left(\gamma \nabla u\mid _{\partial \Omega}\right) \cdot n$, with $n$ the normal to the boundary $\partial \Omega$, are equal to the {\it boundary currents}. Thus the Dirichlet-to-Neumann map $v \mapsto j$ maps boundary voltage potentials to boundary currents, and the correspondence to the discrete setting is provided by identifying $v$ with $\psi_B$, $j$ with $J_B$, the vector $(g_1,\cdots,g_M)$ with $\gamma$, and the Dirichlet-to-Neumann map with the boundary map $J_B=\L_S\psi_B$.

The inverse problem studied by Calderon \cite{calderon} concerns the {\it reconstructibility} of the conductance function $\gamma$ from the knowledge of the Dirichlet-to-Neumann map $v \mapsto j$. Surprisingly, under rather general conditions, the conductance function $\gamma$ is indeed uniquely determined by the Dirichlet-to-Neumann map. See e.g. \cite{borcea} for a review of Calderon's problem.

\smallskip

In the discrete version of a linear resistive network Calderon's inverse problem amounts to the question when and how the conductances $g:=(g_1,\cdots,g_M)$ of the resistors in the resistive network are uniquely determined by the knowledge of the map $\L_S$. Here, somewhat similar to the continuous case, it is throughout assumed that the incidence matrix $\D$ of the circuit graph is {\it known}. 

The discrete version of Calderon's problem has been studied in a number of papers; see in particular \cite{curtis1,curtis2, colin, curtis3}.
Contrary to the continuous case, the values of the conductances are in general {\it not} uniquely determined by the boundary map $J_B=\L_S\psi_B$ for {\it arbitrary} circuit graphs. It {\it is} true for electrical circuits with specific topology, such as the rectangular graphs studied in \cite{curtis1,curtis2} and the circular graphs studied in \cite{colin, curtis3}. Key notion in proving the reconstructibility of the conductance parameters $g_1,\cdots,g_M$ from the boundary map $J_B=\L_S\psi_B$ is the following approach, which is directly extending the approach in the continuous case taken in the original paper by Calderon \cite{calderon}; see also \cite{borcea}.

Consider for a given network graph with incidence matrix $\D$ the map $T$ from the vector $g:=(g_1,\cdots,g_M)$ of conductances to $\L_S$, or equivalently the map from $g$ to the quadratic form defined by $\L_S$, i.e.,
\bq
T: g \mapsto Q_g, \quad Q_g(\psi_B):= \psi_B^T\L_S\psi_B
\eq
Reconstructibility of $g$ from the boundary map $J_B=\L_S\psi_B$ is thus equivalent to {\it injectivity} of the map $T$. This implies the following {\it necessary} condition for reconstructibility. Since the $N_B \times N_B$ Laplacian matrix $\L_S$ is symmetric and has column and row sums zero, it follows that the dimension of the set of all Laplacian matrices $\L_S$ is equal to
$\frac{N_B(N_B-1)}{2}$.
Hence a necessary condition for invertibility of the map $g \mapsto \L_S$, and thus for reconstructibility of $g$, is that 
\bq
M \leq \frac{N_B(N_B-1)}{2}
\eq
In order to derive a {\it sufficient} condition let us recall the equality
\bq
Q_g(\psi_B)= \begin{bmatrix} \psi_B^T & \psi_C^T \end{bmatrix} \L \begin{bmatrix} \psi_B \\ \psi_C \end{bmatrix}, \mbox{ with } \psi_C \mbox{ s.t. }J_C=0
\eq
Then the differential of $T$ at $g$ in the direction of a vector $\kappa \in \mR^M$ is easily seen to be given by the quadratic form
\bq
\left(dT(g)(\kappa)\right)(\psi_B) = \begin{bmatrix} \psi_B^T & \psi_C^T \end{bmatrix} \L^{\kappa} \begin{bmatrix} \psi_B \\ \psi_C \end{bmatrix}, \mbox{ with } \psi_C \mbox{ s.t. }J_C=0,
\eq
where $\L^{\kappa}:=\D [\kappa] \D^T$. Hence the differential $dT(g)$ is injective whenever
\bq
\begin{bmatrix} \psi_B^T & \psi_C^T \end{bmatrix} \L^{\kappa} \begin{bmatrix} \psi_B \\ \psi_C \end{bmatrix}=0 \mbox{ for all } \psi_B,\psi_C \mbox{ s.t. }J_C=0
\eq
implies that $\kappa=0$. Indeed, for rectangular resistive networks as defined in \cite{curtis1, curtis2}, this holds; thus showing reconstructibility of $g$. (Reconstructibility for circular planar graphs was proved in \cite{curtis3} using different methods.) The problem of characterizing all graph topologies (i.e., incidence matrices $\D$) for which the reconstructibility property holds seems to be open.

Let us finally mention that instead of the discrete version of the Dirichlet-to-Neumann map, we may also consider the discrete version of the Neumann-to-Dirichlet map, which is implicitly given by
\bq
J_B \mapsto \psi_B, \quad J_B=\L_S\psi_B, \; \mathds{1}^TJ_B=0
\eq
(Note that indeed the inverse of the map $\L_S$ is well-defined on the subspace of all boundary nodal currents $J_B$ satisfying $\mathds{1}^TJ_B=0$.) Reconstructibility of $g$ remains the same as before.

\section{RLC electrical networks}
The previous framework for linear resistive networks can be directly extended to the {\it steady-state behavior} of linear RLC electrical networks. Indeed, the steady-state behavior of a linear capacitor given by $C\dot{V}=I$ in the frequency-domain is given by the complex relation
\bq
I= j \omega C V, \quad I,V \in \mC,
\eq
with $\omega$ the frequency. Similarly, a linear inductor $L\dot{I}=V$ is described in the frequency-domain as
\bq
I= \frac{1}{j \omega L} V, \quad I,V \in \mC
\eq
Thus the Laplacian matrix $\L_{\mC}$ of the steady-state behavior of an RLC-network with incidence matrix $\D$ is given by
\bq
\label{Lapcomplex}
\L_{\mC} = D G_{\mC} \D^T,
\eq
where $G_{\mC}$ is the complex diagonal matrix with diagonal elements determined by the corresponding edges: real conductances in case of resistors, imaginary numbers $j \omega C$ in case of capacitors, and imaginary numbers $\frac{1}{j \omega L}$ in case of inductors. 

Conversely, somewhat similar to Proposition \ref{prop:laplacian}, it can be seen that any symmetric matrix $\L_{\mC}$ with row and column sums zero, and with off-diagonal elements equal to either $-g$, $-j \omega C$, $-\frac{1}{j \omega L}$, where $g,C,L$ are positive constants, corresponds to an RLC circuit. Indeed, if the $(i,j)$-th element of $\L$ is equal to one of these expressions, then the edge between node $i$ and $j$ corresponds to, respectively, a resistor, capacitor, or inductor. 

On the other hand, taking the Schur complement of such a Laplacian matrix in general leads to a matrix of a more general type; see \cite{willemsverriest2} for some illuminating observations. In general, the problem of characterizing Schur complements of a complex Laplacian matrix as in \eqref{Lapcomplex} seems largely open.

\medskip

A particular case, which has received much attention motivated by direct applications in power networks, is the case where all edges are  inductors ({\it inductive transmission lines}); see e.g. \cite{bolognani, simpsonnature, simpson1,simpson2} and the references quoted therein. In the notation of this paper, the {\it power-flow problem} is to prescribe the vector of {\it complex powers} $[\psi_C] \bar{J}_C$ (with $\bar{}$ denoting {\it complex conjugate}) at the connection nodes (now called the {\it load buses}), as well as the {\it real part} of the complex powers $[\psi_B] \bar{J}_B$ at the boundary nodes (now called the {\it generator buses}), together with the {\it angles} of the complex voltage potentials $\psi_B$. The real part of a complex power $\psi_k \bar{J}_k$ is commonly called the {\it active power}, and the imaginary part the {\it reactive power}. 

\section{Conclusions}
Resistive electrical networks constitute a beautiful example of open, interconnected, large-scale systems, giving rise to an elegant classical mathematical theory, still posing open problems and suggesting important extensions.

\end{document}